\documentclass[12pt]{amsart}
\usepackage{amssymb}
\usepackage{amsthm}
\usepackage{amscd}

\usepackage{amsmath}
\usepackage{epic,eepic}
\usepackage{graphicx}
\theoremstyle{plain}
\newtheorem{lemma}{Lemma}[subsection]
\newtheorem{prop}[lemma]{Proposition}
\newtheorem{thm}[lemma]{Theorem}
\newtheorem{cor}[lemma]{Corollary}
\newtheorem{aplemma}{Lemma~A.\hspace{-1.5mm}}
\newtheorem{approp}{Proposition~A.\hspace{-1.5mm}}
\newtheorem{apthm}{Theorem~A.\hspace{-1.5mm}}
\newtheorem{apcor}{Corollary~A.\hspace{-1.5mm}}
\newtheorem{intthm}{Theorem}
\theoremstyle{definition}

\newtheorem{rema}[lemma]{Remark}

\newtheorem{remb}{Remark}

\newtheorem{defi}[lemma]{Definition}
\newtheorem{exa}[lemma]{Example}
\newtheorem{aprem}{Remark~A.\hspace{-1.5mm}}
\newtheorem{apdefi}{Definition~A.\hspace{-1.5mm}}
\newcommand{\bde}{\begin{defi}}
\newcommand{\ede}{\end{defi}\vspace{1mm}}
\newcommand{\ble}{\begin{lemma}}
\newcommand{\ele}{\end{lemma}}
\newcommand{\bpr}{\begin{prop}}
\newcommand{\epr}{\end{prop}}
\newcommand{\bt}{\begin{thm}}
\newcommand{\et}{\end{thm}}
\newcommand{\bco}{\begin{cor}}
\newcommand{\eco}{\end{cor}}
\newcommand{\bre}{\begin{rema}}
\newcommand{\ere}{\end{rema}}
\newcommand{\brea}{\begin{rema}}
\newcommand{\erea}{\end{rema}\vspace{1mm}}
\newcommand{\breb}{\begin{remb}}
\newcommand{\ereb}{\end{remb}\vspace{1mm}}
\newcommand{\bex}{\begin{exa}}
\newcommand{\eex}{\end{exa}}
\newcommand{\bpf}{\begin{proof}}
\newcommand{\epf}{\end{proof}\vspace{1mm}}

\newcommand{\bade}{\begin{apdefi}}
\newcommand{\eade}{\end{apdefi}}
\newcommand{\bale}{\begin{aplemma}}
\newcommand{\eale}{\end{aplemma}}
\newcommand{\bapr}{\begin{approp}}
\newcommand{\eapr}{\end{approp}}
\newcommand{\bat}{\begin{apthm}}
\newcommand{\eat}{\end{apthm}}
\newcommand{\baco}{\begin{apcor}}
\newcommand{\eaco}{\end{apcor}}
\newcommand{\bare}{\begin{aprem}}
\newcommand{\eare}{\end{aprem}}

\newcommand{\be}{\begin{enumerate}}
\newcommand{\ee}{\end{enumerate}}
\newcommand{\bcd}{\[\begin{CD}}
\newcommand{\ecd}{\end{CD}\]}
\newcommand{\bit}{\begin{itemize}}
\newcommand{\eit}{\end{itemize}}
\newcommand{\bq}{\begin{quote}}
\newcommand{\eq}{\end{quote}}
\newcommand{\ba}{\begin{array}}
\newcommand{\ea}{\end{array}}
\newcommand{\mcE}{\mathcal{E}}
\newcommand{\mcF}{\mathcal{F}}
\newcommand{\mcG}{\mathcal{G}}
\newcommand{\mcH}{\mathcal{H}}

\newcommand{\mcL}{\mathcal{L}}
\newcommand{\mcM}{\mathcal{M}}
\newcommand{\mcN}{\mathcal{N}}
\newcommand{\mcO}{\mathcal{O}}

\newcommand{\mcQ}{\mathcal{Q}}
\newcommand{\mcR}{\mathcal{R}}
\newcommand{\mcS}{\mathcal{S}}

\newcommand{\mcU}{\mathcal{U}}
\newcommand{\mcV}{\mathcal{V}}

\newcommand{\mbC}{\mathbb{C}}

\newcommand{\mbF}{\mathbb{F}}

\newcommand{\mbP}{\mathbb{P}}

\newcommand{\mbZ}{\mathbb{Z}}
\newcommand{\mff}{\mathfrak{f}}

\newcommand{\migi}{\rightarrow}

\newcommand{\isom}{\stackrel{\sim}{\migi}}

\newcommand{\migiincl}{\hookrightarrow}

\newcommand{\migisurj}{\twoheadrightarrow}

\newcommand{\mr}{\mathrm}

\title[The Verschiebung for rank two stable bundles]{An upper bound on the generic degree of  the generalized Verschiebung  for rank two  stable bundles}
\author{Yuichiro Hoshi}
\author{Yasuhiro Wakabayashi}

\address{\emph{Yuichiro Hoshi}
 \newline
 \textnormal{Research Institute for Mathematical Sciences, Kyoto University, Kyoto 606-8502, JAPAN.
 }
 \newline
 \textnormal{\texttt{yuichiro@kurims.kyoto-u.ac.jp}}}

 \address{\emph{Yasuhiro Wakabayashi}
 \newline
 \textnormal{Graduate School of Information Science and Technology, Osaka University, Suita, Osaka 565-0871, JAPAN.}
 \newline
  \textnormal{\texttt{wakabayashi@ist.osaka-u.ac.jp}}}

\date{}
\begin{document}
\maketitle

\footnotetext{2020 {\it Mathematical Subject Classification}: Primary 14H60, Secondary 14D20.}
\footnotetext{Key words:  algebraic curve, stable bundle, moduli space, generalized Verschiebung, Quot scheme.}
\begin{abstract}
In the present  paper, we give an upper bound for the generic degree of the generalized Verschiebung between the moduli spaces of rank two stable bundles with trivial determinant.
\end{abstract}
\tableofcontents 
\section*{Introduction}

Let $k$ be an  algebraically closed field of characteristic $p >0$  and $X$  a  smooth   projective curve over $k$ of genus $g>1$. Denote by $X^{(1)}$ the Frobenius twist of $X$ over $k$. Then, pulling-back  stable bundles on $X^{(1)}$ via the relative  Frobenius morphism $F_{X/k} : X \migi X^{(1)}$ induces  the so-called ``{\it generalized Verschiebung}" rational map
\begin{align} \label{e1}  
\mr{Ver}_{X/k}^n : \mr{SU}_{X^{(1)}/k}^n \dashrightarrow \mr{SU}_{X/k}^n 
\end{align}
 between the moduli spaces
   of rank $n >1$ stable bundles  with trivial determinant on $X^{(1)}$ and $X$ respectively;
   this can be regarded as a higher-rank variant of the Verschiebung between Jacobians.

The geometry  of the rational map $\mr{Ver}_{X/k}^n$, i.e., 
the dynamics of stable bundles with respect to Frobenius pull-back,  has been investigated 
for a long time.
 One motivation  is  the relationship with   representations of the fundamental group of a curve in positive characteristic  (cf., e.g.,  ~\cite{BK1}, ~\cite{Las}, ~\cite{Ya});  indeed, it is well-known that  rank $n$ vector bundles fixed by some powers of the Frobenius morphism come from continuous representations of the fundamental group in $\mr{GL}_n (k)$ (cf. ~\cite{LaSt}).
Also, 
the moduli space $\mr{SU}_{X/k}^n$ and the generalized Verschiebung $\mr{Ver}_{X/k}^n$ are interesting in their own right; other studies regarding  these mathematical objects (e.g., the density of Frobenius-periodic bundles, the loci of Frobenius-destabilized bundles, etc.)
 can be found in  ~\cite{DuMe}, ~\cite{JP}, ~\cite{JRXY}, ~\cite{LaPa}, ~\cite{LaPa2}, ~\cite{LP},    ~\cite{O3},  ~\cite{O2}.

The present paper aims  to  address the generic degree  $\mr{deg}(\mr{Ver}_{X/k}^{n})$ of   $\mr{Ver}_{X/k}^n$ for $n=2$.
The case of $(n, g)  =(2, 2)$ has already been investigated considerably.
We know that, for a genus-$2$ curve $X$,  the compactification of $\mr{SU}_{X/k}^2$ by semistable bundles is canonically isomorphic to the $3$-dimensional projective space $\mbP^3_k$ and the boundary locus can be identified with the Kummer surface associated to $X$.
By this description, 
the rational map $\mr{Ver}_{X/k}^2$ can be expressed  by polynomials of degree $p$ (cf. ~\cite[Proposition A.2]{LaPa2}), which are explicitly described in the cases $p=2$ (cf. ~\cite{LaPa}) and $p=3$ (cf. ~\cite{LaPa2}).
Moreover, the scheme-theoretic base locus of $\mr{Ver}_{X/k}^2$ were computed in  ~\cite[Theorem 2]{LP}.
This result enables us to specify the generic degree $\mr{deg}(\mr{Ver}_{X/k}^{2})$ of
$\mr{Ver}_{X/k}^{2}$, i.e., we have
 $\mr{deg}(\mr{Ver}_{X/k}^{2}) = \frac{p^3+2p}{3}$ (cf. ~\cite[Theorem 1.3]{O3}, ~\cite[Corollary]{LP}). 

However, we have not yet reached a comprehensive understanding of $\mr{Ver}_{X/k}^n$   because
not much seems to be known   for 
general $(p, n, g)$.
That is,  the structure of this rational map remains mysterious despite its own importance! As a step towards understanding it, for example, 
 knowing by what value the generic degree $\mr{deg}(\mr{Ver}_{X/k}^{n})$ is bounded above will be useful information in measuring its complexity.
 (As far as the authors know, it seems that the relevant numerical results are still obtained only when $(n, g)=(2, 2)$.)
 The main result of the present paper, i.e., Theorem \ref{T1} below,  concerns 
 this matter and 
  provides an upper bound  of  the generic degree $\mr{deg}(\mr{Ver}_{X/k}^{2})$    for   infinitely many pairs $(p, g)$.

%---------------------------------------------------------------------[begin theorem]----------------------
\begin{intthm}[= Theorem \ref{T3}] \label{T1}
If $p+1 > g >1$ and $p \neq 2$, then the generic degree $\mr{deg}(\mr{Ver}_{X/k}^2)$ of $\mr{Ver}_{X/k}^2$ satisfies  the following inequality:
\begin{align} \label{e2}
\mr{deg}(\mr{Ver}_{X/k}^2) \leq  p^{g-1} \cdot  \sum_{\theta =1}^{2p-1}\frac{1}{\mr{sin}^{2g-2}(\frac{\pi \cdot \theta}{2 \cdot p})} 
\ \left(=  \sum_{\zeta^{2p} =1, \zeta \neq 1}\frac{(-4p\zeta)^{g-1}}{(\zeta-1)^{2g-2}}\right). 
\end{align}
\end{intthm}
%------------------------------------------------------------------------------[end theorem]-------------

We remark here that an essential ingredient in the proof of the above theorem is the correspondence between the generic fiber of $\mr{Ver}_{X/k}^2$ and a certain Quot scheme.
As discussed by K. Joshi et al. (cf. ~\cite{JP} and an unpublished version\footnote{This version is available at:  {\tt https://arxiv.org/pdf/1311.4359v1.pdf}} of ~\cite{J}),
Quot schemes appear naturally in the study of the action of the Frobenius morphism on vector bundles. 
One of  the new ideas  included in the present paper is to establish 
  this correspondence by  a  dimension  estimate resulting from     Brill-Noether theory (cf. the proof of Lemma \ref{L2}).
We moreover use the  generic \'{e}taleness of $\mr{Ver}_{X/k}^2$ for an ordinary $X$ to   lift the Quot scheme to characteristic $0$.
As a result, the required inequality is obtained by  applying  a formula proved by Y. Holla (cf. ~\cite[Theorem 4.2]{H}), i.e., a special case of the Vafa-Intriligator formula.
The idea of applying Holla's formula to counting objects of this type in positive characteristic is due to K. Joshi (cf. ~\cite{J}, ~\cite{JP}).
A similar argument based on Joshi's idea can be found in 
 the proof of the main theorem in  ~\cite{Wak}.

%----------------------------------------------------------------------[begin subsection]-------------
\subsection*{Notation and Conventions}

Throughout the present paper, we fix  an integer  $g >1$, a prime $p>2$, and an algebraically closed field  $k$   of characteristic $p$.

For each scheme $T$, we shall write $(\mcS ch/T)$ for the category of schemes of finite type over $T$.
 For simplicity, if $R$ is a ring, then we shall write $(\mcS ch/R) := (\mcS ch/\mr{Spec}(R))$.

For a smooth projective curve $X$ over a field $K$, 
we shall write $\Omega_{X/K}$ for  the sheaf of $1$-forms on $X$ relative to $K$.
Also,  for  each  integer $d$, denote by 
$\mr{Pic}_{X/K}^d$ the Picard scheme of $X/K$ classifying isomorphism classes of  line bundles on $X$ of degree $d$.

If $T$ is a scheme and $X_1$,    $X_2$ are  $T$-schemes, then we shall  denote  by $\mr{pr}_i$ ($i=1,2$) the $i$-th projection $X_1 \times_T X_2 \migi X_i$.

\section{Definition of the generalized Verschiebung} \label{S1}\vspace{0mm}

%----------------------------------------------------------------------[begin subsection]-------------
\subsection{} \label{SS2}

Let $R$ be an integrally closed domain of finite type over $k$ 
and $X$ a smooth  projective curve over $R$ of genus $g $.
We shall denote by
\begin{align} \label{e3}  \mcS \mcU_{X/R}^2 : (\mcS ch/R)^\mr{op} \migi (\mcS et)
\end{align}
the set-valued contravariant  functor on the category $(\mcS ch/R)$   which, to any $T \in \mr{Ob}(\mcS ch/R)$, assigns the set of {\it equivalence} classes of $T$-flat families of rank $2$ geometrically stable bundles on $X \times_R T$ with trivial determinant.
Here, 
given two $T$-flat families of geometrically stable bundles $\mcF_1$,  $\mcF_2$ on $X \times_R T$, we say that $\mcF_1$ and $\mcF_2$ are {\it equivalent}  if $\mcF_1 \otimes \mr{pr}_2^*(\mcN) \cong \mcF_2$ for some line bundle $\mcN$ on $T$.
According to ~\cite[Theorem 0.2]{La} and ~\cite[Theorem 9.12]{Lang}, 
 there exists a flat quasi-projective $R$-scheme
\begin{align} \label{e4} \mr{SU}_{X/R}^2 
\end{align}
 that  corepresents universally the functor $\mcS \mcU_{X/R}^2$.
 (Note that the flatness asserted in ~\cite{Lang} is still true even when the base field is of  positive characteristic.)
 The fiber $ \mr{SU}_{X/R}^2 \times_R \xi$ 
 over each geometric point $\xi$ of $\mr{Spec}(R)$
  is irreducible,  smooth, and  of  dimension $3g-3$
  (cf., e.g., ~\cite[Lemma A]{MS2}).

%--------------------------------------------------------------------[begin subsection]-------------
\subsection{} \label{SS3}

Next, write $X^{(1)}$ for the Frobenius twist of $X$ over $R$ and $F_{X/R} : X \migi X^{(1)}$ for the relative Frobenius morphism of $X$ over $R$.
Let 
\begin{align} \label{e5} \mcS\mcU_{X^{(1)}/R}^{2, \circledcirc}
\end{align}
be the subfunctor of $\mcS \mcU_{X^{(1)}/R}^2$ classifying geometrically stable bundles whose pull-back under $F_{X/R}$ is  geometrically stable.
Recall (cf. ~\cite[Proposition 2.3.1]{HL}) that
the property of being geometrically stable is an open condition in flat families.
Hence,  there exists an open subscheme 
\begin{align} \label{e6}    \mr{SU}_{X^{(1)}/R}^{2, \circledcirc} 
\end{align}
of $\mr{SU}_{X^{(1)}/R}^2$ that corepresents universally the functor $\mcS \mcU_{X^{(1)}/R}^{2, \circledcirc}$ (cf. ~\cite[Theorem A.6]{O3}).
The assignment $[\mcF] \mapsto [F_{X/R}^*(\mcF)]$ (for each $[\mcF] \in \mcS \mcU_{X^{(1)}/R}^{2, \circledcirc}$) determines a natural transformation
\begin{align} \label{e7}
   \mcV er_{X/R}^{2, \circledcirc} : \mcS \mcU_{X^{(1)}/R}^{2, \circledcirc}  \migi \mcS \mcU_{X/R}^2
   \end{align}
of functors, which induces a dominant morphism
\begin{align} \label{e8}
   \mr{V er}_{X/R}^{2, \circledcirc} : \mr{SU}_{X^{(1)}/R}^{2, \circledcirc}  \migi \mr{SU}_{X/R}^2
   \end{align}
 between flat $R$-schemes of the same relative dimension (cf. ~\cite{MS2},   ~\cite[Theorem A.6]{O3}).
The formation of $\mr{V er}_{X/R}^{2, \circledcirc}$ is compatible, in an evident sense, with 
restriction to each geometric point of $R$.
Thus, the fact mentioned in the final sentence of the previous subsection  implies  that $\mr{V er}_{X/R}^{2, \circledcirc}$ is generically flat. 

Now, let us specialize the situation  to  the case where $R = k$.
Because of  the  facts mentioned above,
one may define the generic degree of $\mr{V er}_{X/k}^{2, \circledcirc}$, which we denote by 
\begin{align} \label{e9}
\mr{deg} (\mr{Ver}_{X/k}^2).
\end{align}
% and this value
% does not depend on the choice of $X$.
 Here, denote by  $\mcM_g$  the moduli stack classifying smooth projective curves over $k$ of genus $g$.
Since  $\mcM_g$   is 
an irreducible 
 DM stack over $k$ (cf. ~\cite[\S\,5]{DM}), 
 it makes sense to speak of a ``{\it general}" curve, i.e., a curve  that
determines a point of $\mcM_g$ that lies outside a certain fixed closed substack not equal to $\mcM_g$ itself.
%By the flatness of  $ \mr{V er}_{X/R}^{2, \circledcirc}$ (with $R$ as above),
Note that the function on $\mcM_g$ given by  $X  \mapsto \mr{deg} (\mr{Ver}_{X/k}^2)$ is lower semicontinuous.
Therefore, in order to give  an upper bound of the value $\mr{deg} (\mr{Ver}_{X/k}^2)$, 
{\it we are always free to replace $X$ with a general curve in $\mcM_g$}.

%%%%%%%%%%%%%%%%%%%%%%%%%%%%--[ begin  section]---%%%%%%
\vspace{5mm}
\section{Relationship with Quot schemes} 

%----------------------------------------------------------------------[begin subsection]-------------
\subsection{} \label{SS4}
We recall the notion of a Quot scheme as follows.
Let $T$ be a noetherian scheme,
$Y$  a smooth projective  curve over $T$ of genus $g$ and $\mcE$ a vector bundle on $Y$. 
For each integers $r$, $d$ with $r \geq 0$, we shall denote by 
\begin{align} \label{e11} \mcQ uot_{\mcE/Y/T}^{r, d} : (\mcS ch/T)^\mr{op} \migi (\mcS et)
\end{align}
the set-valued contravariant  functor on the category $(\mcS ch/T)$ which, to any $f :T' \migi T$, associates the set of isomorphism classes of 
$\mcO_{Y \times_T T'}$-linear  injections 
 $s :  \mcF \migiincl   (\mr{id}_Y\times f)^*(\mcE)$
such that 
\begin{itemize}
\item
 the cokernel $\mr{Coker}(s)$
  is flat over $T'$
 (which, by the fact that $Y/T$ is smooth of relative dimension $1$, implies that $\mcF$ is {\it locally free}), and 
\item
 $\mcF$ is of rank $r$ and degree $d$. 
\end{itemize}
It  is known (cf. ~\cite[Theorem 5.14]{FGA})  that $\mcQ uot_{\mcE/Y/T}^{r, d}$ 
 may be represented by a proper  scheme over $T$.
 By abuse of notation, write $\mcQ uot_{\mcE/Y/T}^{r, d}$ for the scheme that represents the functor $\mcQ uot_{\mcE/Y/T}^{r, d}$.

%\vspace{5mm}
%----------------------------------------------------------------------[begin subsection]-------------
\subsection{} \label{SS5}
Let
$X$ be a smooth projective curve over $k$ of genus $g$ which is general in $\mcM_g$.
In particular, we may  assume that 
 $X$ is {\it ordinary}.
(Recall that the locus of $\mcM_g$ classifying ordinary curves is open and dense.) 
Moreover,  we assume  that $p+1 > g \left(>1\right)$.
(This assumption will be applied in Lemma \ref{L2} and Proposition \ref{PP3}.)
Let us take a geometric generic point $\eta : \mr{Spec}(K) \migi \mr{SU}_{X/k}^2$ (where $K$ denotes an algebraically closed field over $k$) of 
   $\mr{SU}_{X/k}^2$; this point  classifies a rank $2$ stable bundle $\mcE$ on $X_K := X \times_k K$ with $\mr{det}(\mcE) \cong \mcO_{X_K}$.
Write $X^{(1)}_K$ for the Frobenius twist of $X_K$ over 
$K$ and $F : X_K \migi X_K^{(1)}$ for the relative  Frobenius morphism.
If $\mcG$ is an $\mcO_{X_K^{(1)}}$-module and $\mcH$ is an $\mcO_{X_K}$-module, then the adjunction relation ``$F^*(-) \dashv F_{*}(-)$" gives a natural bijection
\begin{align} \label{e10}
\mr{ad} :  \mr{Hom}_{\mcO_{X_K}}(F^*(\mcG), \mcH) \isom \mr{Hom}_{\mcO_{X^{(1)}_K}}(\mcG, F_{*}(\mcH)),
\end{align}
which is functorial with respect to both $\mcG$ and $\mcH$.  

Since $F$ is  finite and faithfully flat of degree $p$, the direct image $F_*(\mcE)$ forms a vector bundle on $X^{(1)}_K$ of rank $2p$.
Now, consider the Quot scheme 
\begin{align} \label{e12}
  \mcQ   :=  \mcQ uot_{F_{*}(\mcE)/X^{(1)}_K/K}^{2, 0}.   
  \end{align}
This $K$-scheme has the closed subscheme  
\begin{align} \label{e14}
\mcQ^{\mr{triv}}
\ \left(\text{resp.,} \  \mcQ^{\mr{triv}, F} \right)
\end{align}
 classifying injections $s : \mcF \migiincl F_{*} (\mcE)$ with $\mr{det}(\mcF) \cong \mcO_{X^{(1)}_K}$ (resp., $F^*(\mr{det}(\mcF)) \cong \mcO_{X_K}$).
In particular,  there exists a natural closed immersion $\mcQ^{\mr{triv}} \migiincl \mcQ^{\mr{triv},F}$.

%----------------------------------------------------------------------[begin subsection]-------------
\subsection{} \label{SS6}

Since $X$ has been assumed to be ordinary,   $\mr{Ver}^2_{X/k}$ is  generically \'{e}tale
(cf.  ~\cite[Corollary 2.1.1]{MS2}).
Hence, the fiber product $\mr{SU}_{X^{(1)}/k}^{2,\circledcirc} \times_{\mr{SU}^2_{X/k}, \eta} K$ is isomorphic to the disjoint union of finitely many  copies of $\mr{Spec} (K)$.
Let us take   a $K$-rational point $\widetilde{\eta}$ of $\mr{SU}_{X^{(1)}/k}^{2,\circledcirc} \times_{\mr{SU}^2_{X/k}, \eta} K$;
this point 
  corresponds  to a rank $2$ stable bundle $\mcF$  on $X_K^{(1)}$ with trivial determinant.
  The point  $\widetilde{\eta}$ is, by definition,  mapped to $\eta$ by $\mr{Ver}_{X/k}^2$,  meaning that 
  there exists an isomorphism $t : F^* (\mcF) \isom \mcE$.
The morphism $\mr{ad} (t) : \mcF \migi F_{*}(\mcE)$ is   injective because  the composite
\begin{align} \label{e16}
F^*(\mcF) \xrightarrow{F^*(\mr{ad} (t))} F^*(F_{*} (\mcE)) \xrightarrow{\mr{ad}^{-1} (\mr{id}_{F_{*}(\mcE) })} \mcE
\end{align}  
coincides with  $t$ and $F$ is faithfully  flat.
By the stability of $\mcE$ (which implies that $\mr{End}_{\mcO_{X_K}}(\mcE) = K$), 
the $K$-rational point of $\mcQ^{\mr{triv}}$ classifying $\mr{ad}(t)$ does not depend on the choice of $t$.
Hence,
the assignment $\mcF \mapsto \mr{ad} (t)$ determines a well-defined  $K$-morphism
\begin{align} \label{e17}
\mr{SU}_{X^{(1)}/k}^{2,\circledcirc} \times_{\mr{SU}^2_{X/k}, \eta} K \migi \mcQ^{\mr{triv}}.
\end{align}

%-----------------------------------------------------------------------[begin lemma]------------------
\ble  \label{L2}
(Recall that we have assumed that  $p+1 > g>1$.)
Let  us take  an $\mcO_{X_K}$-linear  injection $s : \mcF \migiincl F_{*} (\mcE)$ such that $\mcF$ is a rank $2$ vector bundle.
We shall write $d := \mr{deg}(\mcF)$.
Then, the following assertions hold:
\begin{itemize}
\item[(i)]
The inequality $d \leq 0$ holds, meaning that the maximal degree of rank $2$ vector bundles embedded into $F_*(\mcE)$ is at most $0$.
(A similar maximality assertion can be found in an unpublished version of  ~\cite[Proposition 11.6]{J}.)
\item[(ii)]
The injection 
  $s$ is classified by  $\mcQ^{\mr{triv}, F}$  if and only if  the morphism  $\mr{ad}^{-1} (s)  : F^* (\mcF) \migi \mcE$ is an isomorphism.
  \end{itemize}
 \ele
%-------------------------------------------------------------------------[begin proof]-------------------
\begin{proof}
Let us prove the first assertion.
Suppose that $d > 0$.
By comparing the respective degrees of $F^* (\mcF)$ and $\mcE$,
the (nonzero) morphism $\mr{ad}^{-1}(s)$ cannot be an isomorphism at the generic point of $X_K$.
Hence, since $\mcE$ is locally free of rank $2$ and $X_K$ is a smooth curve over $K$, 
the subsheaf   $\mr{Im}(\mr{ad}^{-1}(s))$ of $\mcE$  forms a line bundle.
Write 
 $\nabla_\mr{can} : F^*(\mcF) \migi \Omega_{X_K/K} \otimes F^*(\mcF)$ for   the 
 connection on $F^*(\mcF)$ determined uniquely by the condition that the sections of the subsheaf $F^{-1}(\mcF)$ are horizontal (cf. ~\cite[Theorem (5.1)]{K}).

In the following, we shall prove the claim that
 the line subbundle  $\mr{Ker} (\mr{ad}^{-1}(s))$ ($\subseteq F^*(\mcF)$)
   is not closed  under  $\nabla_\mr{can}$. 
Suppose, on the contrary,   that  $\mr{Ker} (\mr{ad}^{-1}(s))$ is closed under $\nabla_\mr{can}$.
Since the restriction of $\nabla_\mr{can}$ to $\mr{Ker} (\mr{ad}^{-1}(s))$ has  vanishing $p$-curvature,  $\mr{Ker} (\mr{ad}^{-1}(s))$ is isomorphic to $F^*(\mcU)$ for some line bundle $\mcU$ 
on $X^{(1)}_K$.
Moreover, the resulting (horizontal) composite   $F^*(\mcU) \isom  \mr{Ker} (\mr{ad}^{-1}(s)) \migiincl F^*(\mcF)$ comes, via pull-back by $F$, from an injection $\mcU \migiincl \mcF$.
The composite $F^* (\mcU) \migiincl F^*(\mcF) \xrightarrow{\mr{ad}^{-1}(s)} \mcE$ is identical to the zero map,
 so the corresponding map $\mcU \migiincl \mcF \stackrel{s}{\migiincl} F_{*} (\mcE)$ (via $\mr{ad}$)
must be the zero map.
This contradicts the injectivity of $s$, and hence,  completes the proof of the claim.

Next, observe that 
$F^* (\mcF)$
may be regarded  as  an extension of $\mr{Im} (\mr{ad}^{-1} (s))$  by $\mr{Im} (\mr{ad}^{-1} (s))^{\vee} \otimes \mr{det}(F^*(\mcF))$;
let us fix an isomorphism $\mr{Im} (\mr{ad}^{-1} (s))^{\vee} \otimes \mr{det}(F^*(\mcF))  \isom \mr{Ker} (\mr{ad}^{-1}(s))$.
By  the claim proved above, the following composite  turns out to be injective:
 \begin{align}  \label{e18}
\mr{Im}(\mr{ad}^{-1}(s))^\vee \otimes \mr{det}(F^*(\mcF))  
& \isom \mr{Ker} (\mr{ad}^{-1}(s)) \\
& \migiincl F^*(\mcF) \notag \\
& \hspace{-4mm} \xrightarrow{\nabla_\mr{can}}  \Omega_{X_K/K} \otimes F^*(\mcF) \notag  \\
& \migisurj  \Omega_{X_K/K} \otimes \mr{Im}(\mr{ad}^{-1}(s))  \notag \\
& \isom \mr{pr}_1^*(\Omega_{X/k}) \otimes  \mr{Im}(\mr{ad}^{-1}(s)). \notag
\end{align}
This composite may be verified to be $\mcO_{X_K}$-linear, and hence,  determines a nonzero global section
\begin{align}  \label{e19}
q \in \Gamma (X_K, \mr{Im}(\mr{ad}^{-1}(s))^{\otimes 2}  \otimes \mr{det}(F^*(\mcF))^\vee   \otimes  \mr{pr}_1^*(\Omega_{X/k})). 
\end{align}
Let us fix 
a line subbundle $\mcN$ of $\mcE$  containing the subsheaf  $\mr{Im}(\mr{ad}^{-1}(s))$.
(Hence, $\mcE$ is obtained as an extension of $\mcN^\vee$ by $\mcN$.)
The section $q$ may be regarded, via   the inclusion 
$\mr{Im}(\mr{ad}^{-1}(s))^{\otimes 2} \migiincl \mcN^{\otimes 2}$, as a (nonzero) global section of $\mcN^{\otimes 2}\otimes  \mr{det}(F^*(\mcF))^\vee \otimes \mr{pr}_1^*(\Omega_{X/k})$.
On the other hand, according to ~\cite[Proposition 3.1]{LN2}, the degree of a line subbundle of $\mcE$    is at most $- \frac{g}{2}$ (resp., $-\frac{g-1}{2}$) if $g$ is even (resp., odd). 
This implies
\begin{align}
\mr{deg}(\mcN^{\otimes 2}\otimes  \mr{det}(F^*(\mcF))^\vee \otimes \mr{pr}_1^*(\Omega_{X/k})) & = 2 \cdot  \mr{deg}(\mcN) - p d + 2g-2 \\
&  \leq 2 \cdot \left(-\frac{g-1}{2} \right) - p d + 2g-2   \notag \\
& < 0,  \notag
\end{align}
where the last inequality follows from the assumptions $p+ 1 > g$ and $d>0$.
However, this is a contradiction because  $q \neq 0$.
Consequently, we have  $d \leq 0$, as desired.

We now prove the second assertion.
The ``if\," part is clear from $\mr{det}(\mcE) \cong \mcO_X$, so it suffices to 
 consider the ``only if" part.
Here, we assume that $s$ is classified by $\mcQ^{\mr{triv}, F}$, but   the vector bundle $\mr{Im} (\mr{ad}^{-1} (s)) \left(\neq \{ 0\} \right)$ is of rank $1$.
Just as in the proof of the first assertion,
we can obtain a  nonzero global section
\begin{align} \label{EEE323d}
q \in \Gamma (X_K, \mr{Im}(\mr{ad}^{-1}(s))^{\otimes 2} \otimes \mr{pr}_1^*(\Omega_{X/k}))
\end{align}
of the line bundle $\mr{Im}(\mr{ad}^{-1}(s))^{\otimes 2} \otimes \mr{pr}_1^*(\Omega_{X/k})$ (where we recall that $\mr{det}(F^*(\mcF)) \cong \mcO_{X_K}$).
Moreover, by taking a  line subbundle $\mcN$ of $\mcE$ containing $\mr{Im}(\mr{ad}^{-1}(s))$, we may regard $q$ as a (nonzero) global section of $\mcN^{\otimes 2} \otimes \mr{pr}_1^*(\Omega_{X/k})$.
By a well-known fact of  Brill-Noether theory (cf., e.g., ~\cite{G},  ~\cite{L}),
the existence of such a section  $q$ implies that the scheme-theoretic image
  of the morphism $\mr{Spec} (K) \migi \mr{Pic}^{2g-2 + 2 \cdot \mr{deg}(\mcN)}_{X/k}$ classifying $\mcN^{\otimes 2}\otimes    \mr{pr}_1^*(\Omega_{X/k})$ is of dimension $\leq 2g-2 + 2 \cdot \mr{deg}(\mcN)$ (because of the assumption that   $X$ is general). 
On the other hand, the morphism $\mr{Pic}_{X/k}^{\mr{deg}(\mcN)} \migi \mr{Pic}_{X/k}^{2g-2 + 2 \cdot \mr{deg}(\mcN)}$ determined by $[\mcM] \mapsto [\mcM^{\otimes 2}\otimes   \Omega_{X/k}]$ is finite.
It follows that the scheme-theoretic image
   of the morphism $\mr{Spec} (K) \migi \mr{Pic}_{X/k}^{\mr{deg}(\mcN)}$  classifying $\mcN$  is  of dimension $\leq 2g-2 + 2 \cdot \mr{deg}(\mcN)$.
However, this contradicts the fact proved  in Lemma \ref{L1} below.
Consequently,  the locally free sheaf $ \mr{Im}(\mr{ad}^{-1}(s))$ must be of rank $2$.
Since  $\mr{deg} (F^*(\mcF))$ $= \mr{deg} (\mcE) =0$,   $\mr{ad}^{-1}(s)$ turns out to be  an isomorphism.
This completes the proof of the ``only if\," part, as desired.
\end{proof}
%-----------------------------------------------------------------------[end lemma]-------------------

The following lemma was applied  in the proof of the previous  assertion.

%-----------------------------------------------------------------------[begin lemma]------------------
\ble \label{L1} 
Let us keep the notation in the proof of the second assertion of Lemma \ref{L2}.
Then, the scheme-theoretic image $I$   of the  morphism $\mr{Spec} (K) \migi \mr{Pic}_{X/k}^{\mr{deg}(\mcN)}$ classifying  $\mcN$ is of dimension 
$\geq 2g-1 + 2 \cdot \mr{deg}(\mcN)$. 
 \ele
%-------------------------------------------------------------------------[begin proof]-------------------
\begin{proof}
Define  $D$ to be the moduli space classifying
 isomorphism classes of short exact sequences:
\begin{align} \label{e22}
  \mff_0:  0 \migi \mcN_0 \migi \mcE_0 \migi \mcN^\vee_0 \migi 0
  \end{align}
with  $[\mcN_0] \in I$ and $[\mcE_0] \in \mr{SU}_{X/k}^{2}$.
Let us take
an arbitrary exact sequence $\mff_0$ as above.
By the definition of  stability and the properness of $X/k$, 
 every nonzero endomorphism of  $\mcE_0$ is an isomorphism, and hence, 
\begin{align} \label{e21} h^0 ((\mcN_0^{\vee})^{\vee}\otimes \mcN_0) = h^0((\mcE_0/\mcN_0)^\vee \otimes \mcN_0) = 0.
\end{align}
Thus, by the Riemann-Roch theorem, we obtain 
\begin{align} \label{e20}
 h^1((\mcN^\vee_0)^\vee \otimes \mcN_0) =  g-1 - 2 \cdot \mr{deg}(\mcN_0) = g-1 - 2 \cdot \mr{deg}(\mcN).
 \end{align}
This implies that any fiber of the projection $D \migi I$ given by  $[\mff_0]\mapsto [\mcN_0]$ is  of dimension $\leq \left(g-1 - 2 \cdot \mr{deg}(\mcN)\right)-1 = g-2 - 2 \cdot \mr{deg}(\mcN)$.
Hence, if $\mr{dim}(I) < 2g-1 + 2 \cdot \mr{deg}(\mcN)$,  we have 
\begin{align} \label{e23}\mr{dim}(D)  &  \leq \mr{dim}(I) + 
g-2 - 2 \cdot \mr{deg}(\mcN) \\
 & < (2g-1 + 2 \cdot \mr{deg}(\mcN)) + g-2 - 2 \cdot \mr{deg}(\mcN) 
 \notag \\
& = 3g-3  \left(= \mr{dim}(\mr{SU}_{X/k}^2)\right). \notag
\end{align}
This is a contradiction because
the existence of the extension $0 \migi \mcN \migi \mcE \migi \mcN^\vee \migi 0$ implies that the morphism $D \migi \mr{SU}_{X/k}^2$ given  by $[\mff_0] \mapsto [\mcE_0]$ must be  dominant.
Thus,
we obtain the inequality  $\mr{dim}(I) \geq 2g-1 + 2 \cdot \mr{deg}(\mcN)$, as desired.
\end{proof}
%-----------------------------------------------------------------------[end lemma]-------------------

%-----------------------------------------------------------------------
\begin{rema} \label{R2322}
By an  argument  entirely similar to the proof of the second assertion of Lemma \ref{L2},
we can verify the following assertion:
if $s : \mcF \migiincl F_*(\mcE)$ is an injection classified by a $K$-rational point
of $\mcQ$, then
$s$ is classified by $\mcQ^{\mr{triv}, F}$ if and only if $\mr{det}(\mcF)$ descends  to 
a line bundle on $X$.
\end{rema}
%-----------------------------------------------------------------------

%-----------------------------------------------------------------------
\begin{rema} \label{R2355}
As one of the anonymous  referees  commented,
the discussion after (\ref{EEE323d}) in the proof of 
 %the proof of the ``only if" part in  
    Lemma \ref{L2}, (ii),  would be simplified if the fact    in  ~\cite[Proposition 3.5]{Lau} (see also ~\cite{PaPa}) that a general vector bundle   is  very stable could be applied  (although it is  formulated for  algebraic curves over $\mbC$ in {\it loc.\,cit.}).
In fact, since $\mcE$ is an extension of $\mcN^\vee$ by $\mcN$, 
we can obtain a contradiction as the section $q$  determines a nonzero nilpotent Higgs field on $\mcE$ via the composite injection
\begin{align}
\Gamma (X_K, \mcN^{\otimes 2} \otimes \Omega_{X_K/K}) \migiincl \Gamma (X_K, (\mcE/\mcN)^\vee \otimes \mcE \otimes \Omega_{X_K/K}) \migiincl \mr{Hom}_{\mcO_{X_K}} (\mcE, \mcE \otimes \Omega_{X_K/K}). %\notag
\end{align}
\end{rema}
%-----------------------------------------------------------------------

%--------------------------------------------------[begin proposition]------------------
\bpr \label{P1} 
The morphism (\ref{e17}) is an isomorphism.
In particular, $\mcQ^{\mr{triv}}$ is finite and \'{e}tale  over $K$, and the following equality holds: 
\begin{align} \label{e24}
\mr{deg} (\mr{Ver}^2_{X/k}) = \mr{deg} (\mcQ^{\mr{triv}}/K).
\end{align}
\epr
%-------------------------------------------------------------------------[begin proof]-------------------
\begin{proof}
Let us prove the first assertion.
According to  Lemma \ref{L2} above, 
the morphism (\ref{e17}) induces a bijection between the respective sets of $K$-rational points.
In particular, $\mcQ^{\mr{triv}}$ is finite over $K$.
Let $u : U \migi \mr{Spec}(K)$ be a $K$-scheme defined  as  a connected component of $\mcQ^\mr{triv}$.
Then, the reduced scheme $U_\mr{red}$ associated to $U$ is $\mr{Spec}(K)$.
The $K$-scheme $U$  classifies  an injection 
\begin{align} \label{e25}
s_U : \mcF_U \migi (\mr{id}_{X_K} \times u)^*(F_{*}(\mcE)) \ \left(= (F \times \mr{id}_U)_* ((\mr{id}_{X_K} \times u)^* (\mcE))\right)
\end{align}
on $X_K \times_K U$.
By Lemma \ref{L2} again, the morphism 
\begin{align} \label{e26}
(F \times \mr{id}_U)^* (\mcF_U) \migi (\mr{id}_{X_K} \times u)^* (\mcE)
\end{align}
 corresponding to $s_U$ via the adjunction relation ``$(F\times \mr{id}_U)^*(-) \dashv (F \times \mr{id}_U)_*(-)$"  is surjective  when restricted to $X_K \times_K U_\mr{red}$ ($= X_K$).
By Nakayama's lemma and the fact that both $(F \times \mr{id}_U)^* (\mcF_U)$ and $(\mr{id}_{X_K} \times u)^* (\mcE))$ are locally free,  (\ref{e26}) turns out to be  an isomorphism.
In particular, $\mcF_U$  determines 
 a $K$-morphism
$U \migi \mr{SU}_{X^{(1)}/k}^{2,\circledcirc} \times_{\mr{SU}^2_{X/k}, \eta} K$.
By applying this discussion to  the various connected components of $\mcQ^{\mr{triv}}$,  we obtain a morphism $\mcQ^{\mr{triv}} \migi \mr{SU}_{X^{(1)}/k}^{2,\circledcirc} \times_{\mr{SU}^2_{X/k}, \eta} K$.
One may verify  that this morphism determines, by construction,   the  inverse to (\ref{e17}).
This completes the proof of the first assertion.
The second  assertion follows directly from the first  assertion together with the fact that
$\mr{SU}_{X^{(1)}/k}^{2, \circledcirc} \times_{\mr{SU}_{X/k}^2, \eta}K$ is finite and \'{e}tale over $K$.
\end{proof}
%------------------------------------------[end proposition]-------------------

Also, we obtain the following proposition.

%-----------------------------------[begin proposition]------------------
\bpr \label{P3}
The Quot scheme  $\mcQ$ decomposes into the disjoint union $\mcQ = \mcQ^{\mr{triv}, F} \sqcup \mcR$ for some $K$-scheme $\mcR$.
\epr
%-----------------------------------------------------------------------[begin proof]-------------------
\begin{proof}
Denote by $\mcR$ the closed subscheme of $\mcQ$ classifying injections $s : \mcF \migiincl  F_{*} (\mcE)$ with $\mr{Coker}(\mr{ad}^{-1} (s)) \neq \{0 \}$.
It follows from Lemma \ref{L2} (and its proof) that $\mcQ^{\mr{triv}, F}$ coincides with  the complement of $\mcR$ in $\mcQ$, so it is open in $\mcQ$.
This completes the proof of the assertion.
\end{proof}
%------------------------------------------------[end proposition]-------------------

%----------------------------------------------------------------[begin subsection]-------------
\subsection{} \label{SS7}

Next,  we shall consider the relationship between  $\mcQ^{\mr{triv}}$ and $\mcQ^{\mr{triv}, F}$. 
Denote by
\begin{align} \label{e27}  \mr{Ker}(\mr{Ver}_{X_K/K}^1) 
\end{align}
 the scheme-theoretic inverse image of the point $[\mcO_{X_K}] \in  \mr{Pic}_{X_K/K}^0$  via the classical Verschiebung map 
 $\mr{Ver}_{X_K/K}^1 : \mr{Pic}_{X_K^{(1)}/K}^0 \migi \mr{Pic}_{X_K/K}^0$, i.e., the morphism  given by $[\mcN] \mapsto [F^*(\mcN)]$.
It is well-known
that $\mr{Ker}(\mr{Ver}_{X_K/K}^1)$ is finite and faithfully flat over $K$ of degree $p^g$.
Moreover,  since $X$ has been  assumed to be ordinary, it is  \'{e}tale over $K$, i.e., isomorphic to the disjoint union of $p^g$ copies of $\mr{Spec}(K)$.

%----------------------------------------------------[begin proposition]------------------
\bpr \label{P2} There exists a natural isomorphism
\begin{align} \label{e29}
 \mcQ^{\mr{triv}}  \times_K  \mr{Ker}(\mr{Ver}_{X_K/K}^1) \isom \mcQ^{\mr{triv}, F}
\end{align}
over $K$.
In particular, $\mcQ^{\mr{triv}, F}$ is isomorphic to the disjoint union of  finitely many  copies of  $\mr{Spec}(K)$,  and the following equality holds:
\begin{align} \label{e30}
\mr{deg} ( \mcQ^{\mr{triv}}/K) = \frac{1}{p^g} \cdot \mr{deg} (\mcQ^{\mr{triv}, F}/K).
\end{align}
 \epr
%-----------------------------------------------------------------------[begin proof]-------------------
\begin{proof}
First, we construct 
an action of $\mr{Ker}(\mr{Ver}_{X_K/K}^1)$ on $\mcQ$.
Let us take a pair $([s], [\mcL])$  classified by a  morphism $T \migi \mcQ \times_K  \mr{Ker}(\mr{Ver}_{X_K/K}^1)$,
where $T$
 denotes a $K$-scheme and $s$ denotes an injection $s: \mcF \migiincl \mr{pr}_1^*(F_{*}(\mcE))  \left(= (F \times \mr{id}_T)_* ( \mr{pr}_1^*(\mcE))\right)$ on $X \times_K T$.
Choose  a representative $\mcL$ of $[\mcL]$.
After possibly tensoring it with a suitable line bundle pulled back from $T$, we may assume that
there is  an isomorphism  $\iota : (F \times \mr{id}_T)^*(\mcL) \isom \mcO_{X \times_K T}$.
Then,  the following  composite injection determines a $T$-rational point
$[s \diamond\mcL]$
 of $\mcQ$:
\begin{align} \label{e31} 
s \diamond\mcL : \mcF \otimes \mcL^{\otimes \frac{p+1}{2}} & \xrightarrow{s \otimes \mr{id}}  (F \times \mr{id}_T)_* ( \mr{pr}_1^*(\mcE)) \otimes \mcL^{\otimes \frac{p+1}{2}} \\
& \hspace{4mm}\isom  (F \times \mr{id}_T)_* ( \mr{pr}_1^*(\mcE)\otimes(F \times \mr{id}_T)^*(\mcL^{\otimes \frac{p+1}{2}})) \notag \\
& \hspace{4mm}\isom  (F \times \mr{id}_T)_* ( \mr{pr}_1^*(\mcE)\otimes(F \times \mr{id}_T)^*(\mcL)^{\otimes \frac{p+1}{2}}) \notag \\
& \hspace{4mm}\isom  (F \times \mr{id}_T)_* ( \mr{pr}_1^*(\mcE)), \notag
\end{align}
where the second and last arrows
 are the isomorphisms induced by the projection formula and $\iota$ respectively.
Note that this $T$-rational point does not depend on the choices of the representative $\mcL$ and the isomorphism $\iota$. 
The resulting assignment $([s], [\mcL]) \mapsto [s \diamond\mcL]$  is functorial with respect to $T$, and hence, defines a well-defined action
\begin{align} \label{e32}  \mcQ \times  \mr{Ker}(\mr{Ver}_{X_K/K}^1) \migi \mcQ. \end{align}
Note that if $\mcF$ and $\mcL$ are as above, then
\begin{align} \label{e33}
\mr{det} (\mcF \otimes \mcL^{\otimes \frac{p+1}{2}}) \cong \mr{det} (\mcF) \otimes \mcL^{\otimes 2 \cdot \frac{p+1}{2}} \cong \mr{det} (\mcF) \otimes \mcL^{\otimes (p+1)} \cong \mr{det} (\mcF) \otimes \mcL.
\end{align}
Hence, the action (\ref{e32}) restricts to a morphism
\begin{align} \label{e34}
 \mcQ^{\mr{triv}}  \times_K  \mr{Ker}(\mr{Ver}_{X_K/K}^1) \migi \mcQ^{\mr{triv}, F}.
\end{align}
On the other hand, 
(\ref{e33}) also implies that the assignment
$[s] \mapsto ([s\,\diamond \,\mr{det} (\mcF)^\vee],  [\mr{det} (\mcF)])$  determines the  inverse to (\ref{e34}).
This completes the proof of this proposition.
\end{proof}
%------------------------------------------------------[end proposition]-------------------

%%%%%%%%%%%%%%%%%%%%%%%---[ begin  section]---%%%%%%
\vspace{5mm}
\section{Computation via the Vafa-Intriligator formula} \label{S33} \vspace{3mm}

In this final section, we prove the main theorem by applying an argument similar to the arguments discussed in ~\cite{J} and ~\cite{JP}.
By combining  Propositions \ref{P1} and  \ref{P2},
we obtain the following equalities:
\begin{align} \label{e35}
  \mr{deg}(\mr{Ver}_{X/k}^2) =  \mr{deg}(\mcQ^{\mr{triv}}/K) = \frac{1}{p^g} \cdot \mr{deg}(\mcQ^{\mr{triv}, F}/K).
  \end{align}
Therefore,  to give a bound of    $ \mr{deg}(\mr{Ver}_{X/k}^2)$, it suffices to estimate  the value $\mr{deg}(\mcQ^{\mr{triv}, F}/K)$.

%----------------------------------------------------------------------[begin subsection]-------------
\subsection{} \label{SS9}

In this subsection, we review a numerical formula concerning  the degree of a certain Quot scheme over the field of complex numbers $\mbC$.
Let $C$ be a smooth projective curve over $\mbC$ of genus $g$.
Let $r$, $n$, and $d$ be integers with $1 \leq r \leq n$, and let 
$\mcG$ be 
 a vector bundle  on $C$ of rank $n$ and degree $d$.
Then,  we define invariants
 \begin{align} \label{e36}e_{\mr{max}}(\mcG,r) & := \mr{max} \big\{  \mr{deg}(\mcF) \in \mbZ\, \big|\, \text{$\mcF$ is a subbundle of $\mcG$ of rank $r$ }   \big\},   \\
  s_r(\mcG) & := d \cdot r - n \cdot e_{\mr{max}}(\mcG, r). \notag \end{align}
Denote by $\mr{U}_C^{n,d}$ the moduli space of stable bundles on $C$ of rank $n$ and degree $d$.
Since  $\mr{U}_C^{n,d}$ is irreducible (cf., e.g., ~\cite[Lemma A]{MS2}),
it makes sense to speak of a ``general" stable bundle in $\mr{U}_C^{n,d}$, i.e., a stable bundle that corresponds to a point of the scheme $\mr{U}_C^{n,d}$ that lies outside a certain  fixed closed subscheme.
If $\mcG$ is  a general stable bundle in $\mr{U}_C^{n,d}$,
then it holds (cf. ~\cite{HIR}, ~\cite[\S\,1]{LN}) that
$ s_r (\mcG) = r (n-r)(g-1) + \epsilon$, where
$\epsilon$ is the 
integer uniquely determined by the equality just before and 
$ 0 \leq \epsilon < n$.
Also, 
the number $\epsilon$
coincides (cf. ~\cite[\S\,1]{H}) with the dimension of every irreducible component of the Quot scheme 
$\mcQ uot^{r, e_\mr{max}(\mcG,r)}_{\mcG/C/\mbC}$.
If, moreover, the equality $s_r(\mcG) =r(n-r)(g-1)$ holds (i.e., $\mr{dim}(\mcQ uot^{r, e_\mr{max}(\mcG, r)}_{\mcG/C/\mbC})=0$), then the Quot scheme $\mcQ uot^{r, e_\mr{max}(\mcG,r)}_{\mcG/C/\mbC}$ for $\mcG$ general is \'{e}tale over $\mbC$   (cf. ~\cite[Proposition 4.1]{H}).
Finally,  under this particular assumption, a formula for the degree of this Quot scheme was given by Holla
as follows. 
%-----------------------------------------------------------------------[begin theorem]------------------
\bt   \label{T2}
  Let $C$ be a smooth projective curve over $\mbC$ of genus $g$ and $\mcG$  a general stable bundle in $\mr{U}_C^{n,d}$.
Write $(a,b)$ for the unique pair of integers such that $d =an-b$ with $0 \leq b < n$.
Also, we suppose that the equality $s_r(\mcG) = r(n-r)(g-1)$ (equivalently, $e_\mr{max}(\mcG,r) = (dr - r(n-r)(g-1))/n$) holds.
 Then, the degree  $\mr{deg}(\mcQ uot^{r,e_\mr{max}(\mcG,r)}_{\mcG/C/\mbC}/\mbC) $ of $\mcQ uot^{r,e_\mr{max}(\mcG,r)}_{\mcG/C/\mbC}$ over $\mbC$ is calculated by the following formula.
 \begin{align} \label{e37}  \frac{(-1)^{(r-1)(br-(g-1)r^2)/n}n^{r(g-1)}}{r!} \cdot \sum_{\zeta_1 , \cdots , \zeta_r}\frac{\prod_{i=1}^r \zeta_i^{b-g+1}}{\prod_{i \neq j} (\zeta_i - \zeta_j)^{g-1}},
 \end{align}
 where 
 the sum  is taken over the set of $r$-tuples $(\zeta_1, \cdots, \zeta_r) \in \mbC^r$ of mutually distinct $n$-th roots of unity in $\mbC$.
 \et
%-----------------------------------------------------------------------[begin proof]-------------------
\begin{proof}
The assertion follows from  ~\cite[Theorem 4.2]{H}, where the ``$k$'' (resp., ``$r$'') corresponds to our $r$ (resp., $n$).
\end{proof}
%-----------------------------------------------------------------------[end theorem]-------------------
 
%----------------------------------------------------------------------[begin subsection]-------------
\subsection{} \label{SS8}

With  the notation in the previous section,
 we  relate the above formula  
 to the degree of the related  Quot schemes, and then, 
give
an upper bound of the value  $\mr{deg} (\mcQ^{\mr{triv}, F}/K)$.

%----------------------------------------------------[begin proposition]------------------
\bpr \label{PP3} 
(Recall that we have assumed that $p+1 > g>1$.)
We have the following inequality:
\begin{align} \label{e38}
\mr{deg} (\mcQ^{\mr{triv}, F}/K) \leq  p^{2g-1} \cdot  \sum_{\theta =1}^{2p-1}\frac{1}{\mr{sin}^{2g-2}
\left(\frac{\pi \cdot  \theta}{2 \cdot p}\right)}. 
\end{align}
 \epr
%---------------------------------------------------[begin proof]-------------------
\begin{proof}
By Propositions \ref{P1}, \ref{P2}, and the italicized comment at the end of \S\,\ref{S1}, 
 one  may assume, without loss of generality,  that
 $k$ is an algebraic closure of $\mbF_p := \mbZ/p \mbZ$ and the transcendence  degree of the  field extension $K/k$ is finite.
Denote by $W$ the ring of Witt vectors with coefficients in $K$ and
$L$ the fraction field of $W$.
The above assumption implies that we can find an injective morphism of fields $L \migiincl  \mbC$.
Since 
$\mr{dim}(X^{(1)}_K) =1$, which implies  $H^2(X_K^{(1)}, \Omega^\vee_{X_K^{(1)}/K})=0$,
 it follows from  well-known generalities on  deformation theory  that 
  $X^{(1)}_K$ may be lifted to a smooth projective  curve $X^{(1)}_W$ over $W$ of genus $g$.
In a similar vein, the equality  $H^2(X^{(1)}_K, \mcE nd_{\mcO_{X_K^{(1)}}} (F_{*}(\mcE)))=0$ implies that  $F_{*}(\mcE)$ may be lifted to a vector bundle $\mcV_W$ on $X^{(1)}_W$.
Now let  $v$ be a $K$-rational point of  $\mcQ^{\mr{triv}, F}$, which  classifies 
an injection  $s : \mcF \migiincl F_{*}(\mcE)$.
By Proposition \ref{P3}, the tangent space of $\mcQ^{\mr{triv}, F}$ at $v$ may be  identified with the tangent space of $\mcQ$ at the same point, so it is isomorphic to 
the $K$-vector space $\mr{Hom}_{\mcO_{X_K^{(1)}}}(\mcF, \mr{Coker} (s))$.
Also,  the obstruction to lifting $v$ to any first order thickening of $\mr{Spec}(K)$ is given by an element of $\mr{Ext}^1_{\mcO_{X_K^{(1)}}}(\mcF,  \mr{Coker}(s))$.
On the other hand,  the \'{e}taleness of  $\mcQ^{\mr{triv}, F}/K$   (cf. Proposition \ref{P2}) implies 
 the equality  $\mr{Hom}_{\mcO_{X^{(1)}}}(\mcF, \mr{Coker}(s))=0$,  and hence, we have   $\mr{Ext}^1_{\mcO_{X^{(1)}}}(\mcF, \mr{Coker}(s))=0$ by Lemma \ref{L4}  below.
Thus, it follows that
 $v$ may be  lifted uniquely  to a $W$-rational point of $\mcQ uot_{\mcV_W/X^{(1)}_W/W}^{2, 0}$.
 In particular, there exists an open and closed subscheme $\mcQ_W^{\mr{triv}, F}$ of $\mcQ uot_{\mcV_W/X^{(1)}_W/W}^{2, 0}$  whose special fiber coincides with  $\mcQ^{\mr{triv}, F}$.
% Here,  it follows from  a routine argument 
% that  $L$ may be supposed to be a subfield of $\mbC$.
Write  $X^{(1)}_{\mbC}$ for   the base-change of $X^{(1)}_{W}$ via the morphism $\mr{Spec}(\mbC) \migi  \mr{Spec}(W)$ induced by the composite embedding $W \migiincl L \migiincl  \mbC$, and 
$\mcV_\mbC$ for the pull-back of $\mcV_W$ via the natural morphism $X^{(1)}_{\mbC} \migi X^{(1)}_{W}$. 
Then the  degree of $\mcV_\mbC$ coincides with the degree of   $F_*(\mcE)$, 
 so $\mcV_\mbC$ is a  vector bundle of degree $\delta  := 2 \cdot (p-1)(g-1)$ (cf. the proof of Lemma \ref{L4} below). 
Since
 $F_{*}(\mcE)$ is stable (cf. \cite[Theorem 2.2]{SUN}),
one may  verify   from the definition of stability and the properness of Quot schemes (cf.  ~\cite[Theorem 5.14]{FGA}) that $\mcV_\mbC$ is a stable vector bundle.
By the first assertion of  Lemma \ref{L2}, together with the properness of  $\mcQ_W^{\mr{triv}, F}/W$,
 $\mcQ uot_{\mcV_\mbC/X^{(1)}_\mbC/\mbC}^{2, 0}$ classifies maximal subbundles of $\mcV_\mbC$.
 One may assume, without loss of generality, that  the deformation $\mcV_\mbC$ is   sufficiently general  in
$\mr{U}_{X_\mbC^{(1)}}^{2p, \delta}$  so that 
 the dimension of  any component in  $\mcQ uot^{2, 0}_{\mcV_\mbC/X^{(1)}_\mbC/\mbC}$ is the same (cf. \S\,\ref{SS9}).
 (In fact, let $\overline{\mr{U}}_{X_K^{(1)}}^{2p, \delta}$ and $\overline{\mr{U}}_{X_W^{(1)}}^{2p, \delta}$ denote the moduli spaces of rank $2p$ and degree $\delta$ semistable bundles on $X^{(1)}_K/K$ and  $X_W^{(1)}/W$, respectively.
 Then, by applying, e.g., ~\cite[Lemma A]{MS2}, one may verify that the inverse image via the specialization map $\overline{\mr{U}}_{X_W^{(1)}}^{2p, \delta} (L) \migi \overline{\mr{U}}_{X_K^{(1)}}^{2p, \delta} (K)$ of the point classifying $F_*(\mcE)$ is dense in $\mr{U}_{X_\mbC^{(1)}}^{2p, \delta} (\mbC) \subseteq \overline{\mr{U}}_{X_W^{(1)}}^{2p, \delta} (\mbC)$.
 For the moduli space of semistable bundles in mixed characteristic, we refer to ~\cite{La2}.)
In particular, if we  write $\mcQ_\mbC^{\mr{triv}, F} := \mcQ_W^{\mr{triv}, F} \times_W \mbC$,
then the finiteness of  $\mcQ^{\mr{triv}, F}/K$  implies that  $\mcQ_\mbC^{\mr{triv}, F}$,  hence also  $\mcQ uot_{\mcV_\mbC/X^{(1)}_\mbC/\mbC}^{2, 0}$,  is $0$-dimensional.
  Thus, we have
\begin{align} \label{e44}
\mr{deg} (\mcQ^{\mr{triv}, F}/K)  
 = 
\mr{deg} (\mcQ_W^{\mr{triv}, F}/W)  
 =   \mr{deg} (\mcQ_\mbC^{\mr{triv}, F}/\mbC) \leq \mr{deg} (\mcQ uot^{2, 0}_{\mcV_\mbC/X^{(1)}_\mbC/\mbC}/\mbC).
\end{align}
If, moreover,  $(a,b)$ is  the unique pair of integers  satisfying  $\mr{deg}(\mcV_\mbC ) = 2 p \cdot a -b$ with $0 \leq b < 2p$,
 then it follows from the hypothesis $p+1>g$ that $a = g-1$ and $b= 2 (g-1)$.
Thus, since $\mcV_\mbC$ is assumed to be general, 
we can  apply Theorem \ref{T2} in the case where the data 
``$(C, \mcG, n, d, r, a, b, e_{\mr{max}}(\mcG, r))$"
 is taken to be
 \begin{align} \label{e42}
 (X^{(1)}_{\mbC}, \mcV_\mbC, 2p, \delta, 2, g-1, 2(g-1), 0)
 \end{align}
  and   obtain the following sequence of equalities
\begin{align}  \label{e43}
 & \mr{deg}_\mbC(\mcQ uot^{2, 0}_{\mcV_\mbC/X^{(1)}_\mbC/\mbC})  \\
 = & \ \frac{(-1)^{(2-1)(2(g-1)2-(g-1)2^2)/2p}(2p)^{2(g-1)}}{2!} \cdot \sum_{\rho_1, \rho_2}\frac{\prod_{i=1}^2 \rho_i^{2(g-1)-g+1}}{\prod_{i \neq j} (\rho_i - \rho_j)^{g-1}} \notag \\
= & \ (-1)^{g-1}\cdot 2^{2g-2} \cdot p^{2g-1} \cdot  \sum_{\zeta^{2p}=1, \zeta \neq 1}\frac{\zeta^{g-1}}{(\zeta-1)^{2g-2}} \notag  \\
= & \ 2^{g-1} \cdot p^{2g-1} \cdot  \sum_{\zeta^{2p}=1, \zeta \neq 1}\frac{1}{(1-\frac{\zeta + \zeta^{-1 }}{2})^{g-1}}  \notag \\
= & \ p^{2g-1} \cdot  \sum_{\theta =1}^{2p-1}\frac{1}{\mr{sin}^{2g-2}(\frac{\pi \cdot \theta}{2 \cdot p})}.\notag
\end{align} 
Thus, the assertion follows from (\ref{e44}) and (\ref{e43}).
\end{proof}
%-----------------------------------------------------------------------[end proposition]-------------------
The following lemma was applied  in the proof of the previous proposition.

 %-----------------------------------------------------------------------[begin lemma]------------------
\ble \label{L4}  Let $s:\mcF \migiincl F_{*}(\mcE)$ be  an injection classified by  a $K$-rational point of $\mcQ^{\mr{triv}, F}$.
 Write $\mcG := \mr{Coker}(s)$.
 Then 
 $\mcG$ is a vector bundle on $X^{(1)}_K$,  and
  the following equality holds:
 \begin{align} \label{e46}\mr{dim}_K(\mr{Hom}_{\mcO_{X_K^{(1)}}}(\mcF, \mcG) )= \mr{dim}_K(\mr{Ext}^1_{\mcO_{X_K^{(1)}}}(\mcF, \mcG)). \end{align}
 \ele
%-----------------------------------------------------------------------[begin proof]-------------------
\begin{proof}
First, we verify that $\mcG$ is  a vector bundle.
Recall (cf. Lemma \ref{L2}) that 
the composite 
\begin{align} \label{e45}
F^*(\mcF) \xrightarrow{F^*(s)} F^*(F_{*}(\mcE))  \xrightarrow{\mr{ad}^{-1}(\mr{id}_{F_{*}(\mcE)})} \mcE
\end{align}
  is an isomorphism.
Hence,  
the composite $\mr{Ker}(\mr{ad}^{-1}(\mr{id}_{F_{*}(\mcE)})) \migiincl F^*(F_{*}(\mcE)) \migisurj F^*(\mcG)$ is an isomorphism,  so  $F^* (\mcG)$ is a vector bundle.
By the faithful flatness of $F$, $\mcG$ turns out to be a vector bundle on $X_K^{(1)}$, as desired.

Next we shall prove (\ref{e46}).
Since  $F$ is finite,  we have  an equality of Euler characteristics $\chi(F_{*}(\mcE)) = \chi (\mcE) = 2(1-g)$.
Since $ \mr{rk} (\mcH om_{\mcO_{X_K^{(1)}}}(\mcF, \mcG)) = 2 \cdot (2p-2)$, it follows from 
the Riemann-Roch theorem that 
\begin{align} \label{e47} 
\mr{deg}(F_{*}(\mcE)) =
% & \ 
\chi (F_*(\mcE)) - \mr{rk}(F_{*}(\mcE))(1-g) 
%\\
= 
%& \
 \chi (\mcE) - 2  p \cdot (1-g)
 %\notag \\
= 
%& \
 2 \cdot (p-1)(g-1), 
 %\notag
\end{align}
and that 
\begin{align} \label{e48}
 \mr{deg}( \mcH om_{\mcO_{X_K^{(1)}}}(\mcF, \mcG)) = 
 & \
  2 \cdot \mr{deg}(\mcG) - (2p-2) \cdot \mr{deg}(\mcF) 
   \\
  = 
  & \
   2 \cdot \mr{deg}(F_{*}(\mcE)) 
   - 0 
   \notag \\
  =
   &  \
    4 \cdot  (p-1) (g-1).  
    %\hspace{-5mm}
    \notag
  \end{align}
Finally,
by applying the Riemann-Roch theorem again, we obtain 
\begin{align} \label{e49} 
& \  \mr{dim}_K(\mr{Hom}_{\mcO_{X_K^{(1)}}}(\mcF, \mcG)) - \mr{dim}_K(\mr{Ext}_{\mcO_{X_K^{(1)}}}^1(\mcF, \mcG)) \\
 = &  \ \mr{deg}(\mcH om_{\mcO_{X_K^{(1)}}} (\mcF, \mcG)) + \mr{rk} (\mcH om_{\mcO_{X_K^{(1)}}} (\mcF, \mcG)) (1-g) \notag \\
 =  & \ 4 \cdot (p-1)(g-1) + 2\cdot (2p-2)  (1-g) \notag \\
 =  & \ 0, \notag \end{align}
 thus completing the proof of this lemma.
\end{proof}
%-----------------------------------------------------------------------[end lemma]-------------------

By applying the results obtained so far, we  conclude the main result of the present paper.

%-----------------------------------------------------------------------[begin theorem]------------------
\bt [= Theorem \ref{T1}] \label{T3}
 Let $X$ be a smooth projective curve over  $k$  of genus $g$ with   $p+1>g >1$ and $p \neq 2$. Then,
the following inequality holds:
\begin{align} \label{e50}
\mr{deg}(\mr{Ver}_{X/k}^2) \leq  p^{g-1} \cdot  \sum_{\theta =1}^{2p-1}\frac{1}{\mr{sin}^{2g-2}(\frac{\pi \cdot \theta}{2 \cdot p})} \ \left(=  \sum_{\zeta^{2p} =1, \zeta \neq 1}\frac{(-4p\zeta)^{g-1}}{(\zeta-1)^{2g-2}}\right).
\end{align}
\et
%-----------------------------------------------------------------------[begin proof]-------------------
\begin{proof}
By the italicized comment at the end of \S\,\ref{S1}, 
 one  may assume, without loss of generality,  that $X$ is sufficiently  general  for which
 the above discussions work.
Then, by the discussion at the beginning of \S\,\ref{S33} and  Proposition \ref{PP3},  we have
 \begin{align} \label{e51}
  \mr{deg}(\mr{Ver}_{X/k}^2)  =  \frac{1}{p^g} \cdot   \mr{deg}(\mcQ^{\mr{triv}, F}/K) \leq  p^{g-1} \cdot  \sum_{\theta =1}^{2p-1}\frac{1}{\mr{sin}^{2g-2}\left(\frac{\pi \cdot \theta}{2 \cdot p}\right)}.
 \end{align}
This completes the proof of the theorem.
 \end{proof}
%-----------------------------------------------------------------------[end theorem]-------------------

%---------------------------------------------------------------------[begin remark]-------------------
% \begin{rema}
 Similarly to the discussion in  ~\cite[\S\,6.2, (2)]{Wak}, we can describe the right-hand side of (\ref{e50}) as a polynomial with respect to $p$ of degree $3g-3$.
 For example,  (\ref{e50}) reads 
 \begin{align}
 \mr{deg}(\mr{Ver}_{X/k}^2) \leq \frac{4p^3-p}{3} \ \text{if $g=2$,} \hspace{3mm} \text{and}
 \hspace{3mm} \mr{deg}(\mr{Ver}_{X/k}^2) \leq \frac{16p^6 +40 p^4 -11p^2}{45} \ \text{if $g=3$}.
 \end{align}
 By comparing with  the explicit computation of $\mr{deg}(\mr{Ver}_{X/k}^2)$ for $g=2$ obtained already   (cf. Introduction), we see that (\ref{e50}) is not optimal.
 In particular,  (by considering the discussion in the proof of Proposition \ref{PP3}) 
 the $K$-scheme $\mcR \left(= \mcQ \setminus   \mcQ^{\mr{triv}, F}\right)$  in Proposition \ref{P3}  for any general  genus-$2$ curve  turns out to be nonempty.
 Hence, the comment in Remark \ref{R2322} implies the following assertion.
% \end{rema}
 %---------------------------------------------------------------------[end remark]-------------------
 
 %-----------------------------------------------------------------------[begin theorem]------------------
\bco  \label{Cffhg2}
 Let $X$ be a smooth projective curve over  $k$  of genus $2$ and $\eta : \mr{Spec}(K) \migi \mr{SU}_{X/k}^2$ (where $K$ denotes an algebraically closed field over $k$)
 a geometric generic point of $\mr{SU}_{X/k}^2$.
 Denote by $\mcE$ the rank $2$ stable bundle on $X_K := X \times_k K$ classified by $\eta$ and by $F$ the relative Frobenius morphism $X_K \migi X_K^{(1)}$ of $X_K$ over $K$.
 Suppose that $X$ is general in $\mcM_2$.
Then,  the direct image $F_* (\mcE)$ of $\mcE$ via $F$ admits an $\mcO_{X_K^{(1)}}$-submodule $\mcF$ which  is a rank $2$ vector bundle of degree $0$  and whose determinant $\mr{det}(\mcF)$ does not descend to a line bundle on $X$.
\eco
%-----------------------------------------------------------------------[begin proof]-------------------
 
 %--------------------------------------------[begin subsection]-------------
\subsection*{Declaration of competing interest}
The authors declare that they have no known competing financial interests or personal relationships that could have appeared to influence the work reported in this paper.
  
%--------------------------------------------[begin subsection]-------------
\subsection*{Acknowledgements}
We would like to thank Professors Kirti Joshi, Akio Tamagawa, and the  anonymous referees 
for their productive comments and suggestions  which helped improving the quality of the paper.
Also, we are grateful for the many constructive conversations we had with the moduli space of  rank $2$ stable bundles $\mr{SU}_{X/k}^2$, who live in the world of mathematics!
The first author was partially supported by JSPS KAKENHI Grant Number 21K03162.
The second author was partially supported by 
 JSPS KAKENHI Grant Numbers
18K13385, 21K13770.
This work was supported by the Research Institute for Mathematical Sciences, an International Joint Usage/Research Center located in Kyoto University.

%%%%%%%%%%%%%%%%%%%%%%%%%%%%%%%%%%%%%%%%%
\end{document}